\documentclass[10pt,a4paper,final]{article}
\usepackage[utf8]{inputenc}
\usepackage[english]{babel}
\usepackage{amsmath}
\usepackage{amsfonts}
\usepackage{amssymb}
\usepackage{graphicx}
\usepackage{multirow}
\usepackage{revstat}
\usepackage{subcaption}
\newcommand{\ud}{\mathrm{d}}
\begin{document}
\title{Information Complexity Criterion for Model Selection in Robust Regression Using A New Robust Penalty Term}
\renewcommand{\titleheading}
             {}  
\author{\authoraddress{Esra Pamuk\c{c}u}
                      {(1) F{\i}rat University, Faculty of Science, Department of Statistics, Elaz{\i}\v{g}, Turkey,                       \\ (epamukcu@firat.edu.tr)}
                       \\
                        \authoraddress{Mehmet Niyazi \c{C}ankaya}
                      {Department of International Trading and Finance, Faculty of Applied Sciences, U\c{s}ak University \&
                      	\\  Department of Statistics, Faculty of Arts and Science,
                      	U\c{s}ak University,  U\c{s}ak 
                       \\ (mehmet.cankaya@usak.edu.tr *Corresponding Author)}
      }

\maketitle

\begin{abstract}
Model selection is basically a process of finding the best model from the subset of models in which the explanatory variables are effective on the response variable. The log likelihood function for the lack of fit term and a specified penalty term are used as two parts in a model selection criteria. In this paper, we derive a new tool for the model selection in robust regression. We introduce a new definition of relative entropy based on objective functions. Due to the analytical simplicity, we use Huber's objective function $\rho_H$ and propose our specified penalty term $C_0^{\rho_H}$ to derive new Information Complexity Criterion ($RICOMP_{C_0^{\rho_H}}$) as a robust model selection tool. Additionally, by using the properties of $C_0^{\rho_H}$, we propose a new value  of tuning parameter  called $k_{C_0}$ for the Huber's  $\rho_H$. If a contamination to normal distribution exists, $RICOMP_{C_0^{\rho_H}}$ chooses the true model better than the rival ones.  Monte Carlo Simulation studies are carried out to show the utility both of  $k_{C_0}$ and  $RICOMP_{C_0^{\rho_H}}$. A real data example  is also given. 
\end{abstract}

\begin{keywords}
Information theory, Model selection, Robust statistics.
\end{keywords}



\section{Introduction}  
\label{Sec1}              

Since the early 1970s, it has been possible to come across many studies on the model selection algorithm and information criteria. These include classical methods and methods based on information criteria. Classical model selection methods are generally performed by hypothesis testing. However, most statisticians and some other scientists stress that the confidence interval which is used in the classical model selection methods is baseless \cite{Linzuc86}. Some scientists have stated that the hypothesis testing approach does not have a theoretical accuracy and is generally insufficiently valid \cite{Bhandown1977}.

Model selection procedures based on information criteria (IC) provide an alternative to classical approaches. The basic idea of such criteria is to minimize the Kullback-Leibler distance between the true model distribution and the distribution of response variable which belongs to most appropriate model. \cite{Akaike73} and \cite{Akaike74} with the successive published articles  has led to developments in the field of statistical modeling and statistical model selection or evaluation. Therefore, he is considered as a first researcher in this field. Akaike Information Criterion (AIC) proposed by Akaike is given by
\begin{equation}\label{AIC}
AIC=-2\log L(v;\hat{\boldsymbol{\beta}},\hat{\sigma}) + 2p,
\end{equation}
\noindent where $\hat{\boldsymbol{\beta}} $ is a vector of parameters $\beta_1,\beta_2,\cdots,\beta_p$ and $ \sigma $ is scale parameter. $\hat{\boldsymbol{\beta}}$ and $\hat{\sigma}$ are ML estimators of these parameters of probability density (p.d.) function $f(v;\boldsymbol{\beta},\sigma)$. $p$ is the number of estimable parameters.

In the following years, many model selection criteria based on AIC were proposed. Some of these are Takeuchi's Information Criteria also known as the Generalized Akaike Information Criteria  proposed by \cite{takeuchi1976}, Corrected Information Criteria  developed by \cite{sugiura1978}, the Bayesian Information Criterion also known as the Bayesian Information Criterion developed by \cite{Schwarz78}, Hannan-Quinn Information Criterion developed by \cite{HanQuin1979}, Consistent Akaike Information Criterion developed by \cite{Bozdogan87}.

ICOMP (I; Information-COMP; Complexity) which are information criteria developed by \cite{Bozdogan1988} is a tool for model selection in multivariate linear and nonlinear models. Whereas AIC is only intended to strike a balance between the lack of fit and the penalty term, ICOMP aims to establish this balance by taking into account a complexity measure that measures how the parameters in the model relate to each others. Therefore, instead of penalizing the number of independent parameters directly, it penalizes the covariance complexity of the model. In the most general form, ICOMP is given by
\begin{equation}\label{AICcp}
ICOMP=-2\log L(v;\hat{\boldsymbol{\beta}},\hat{\sigma}) + C(\Sigma).
\end{equation}
The lack of fit part $-2\log L$ remains unchanged while the penalty term $C(\Sigma)$ can be chosen by user with different ways  to increase the performance of ICOMP. 

ICOMP is defined as follow, 
\begin{equation}\label{ICOMPC1}
ICOMP=-2\log L(v;\hat{\boldsymbol{\beta}},\hat{\sigma}) + 2C_1(\Sigma(\hat{\boldsymbol{\beta}})),
\end{equation}
\noindent where $\Sigma$ is the estimated covariance matrix of the model and $C_1 (\Sigma)$  is the maximal entropic complexity defined by \cite{Bozdogan1988}. We will provide clear explanations about the complexity $C_1$ (see section \ref{relativeentropyandcomplexity}).

There are several forms and theoretical justifications of ICOMP defined as a general form in Eq.  \eqref{ICOMPC1}. So far, ICOMP have been compared with other model selection criteria on many different methods and it has been observed that they give better results clearly \cite{Kwonetal08,AkbilgicBoz11,Howeetal12,KocBoz15,Pamukcu15}.

In this paper, our main idea is to produce a new robust model selection criterion when a contamination into normal distribution occurs. We follow $C_0$ based on work of \cite{Emden71} and introduce our penalty term $C_0^{\rho_H}$ by using Huber's objective function ($\rho_H$)  in the framework of $C_0$. We choose Huber's $\rho$ function because of its analytical tractability for calculating the expectation $E(\cdot)$. The normal distribution is chosen as the underlying distribution for calculating $E(\cdot)$ in Shannon's relative entropy or Shannon's divergence. 

The rest of paper is organized as follows. Section \ref{robustcriteria} introduces  the robust model selection criteria. Section \ref{relativeentropyandcomplexity} provides the theoretical baseline of the information  complexity   criteria based on Shannon's relative entropy. We propose our penalty term $C_0^{\rho_H}$ in Section \ref{hubercankayapamukcu}. Section \ref{kc0new} presents the proposition of the new value of Huber's tuning parameter $k$. In section \ref{ourICOMP}, we propose the robust ICOMP based on Huber's $\rho_H$ function. The numerical examples are demonstrated in Section \ref{numericalexamples}. Finally, Section \ref{conclusion} is divided by conclusion and brief discussion.

\section{Material and Method}\label{matemet} 
\subsection{Robust Model Selection}\label{robustcriteria}
Over the past 50 years, many robust estimation techniques have been proposed as alternative approaches to the  sum of least squares technique for regression modeling when there are small deviations from the model assumption about underlying distribution and unusual observations in the data sets \cite{Andrewsetalloc72,Huber81,Hampeletal86}.

Whereas there are a lot of works about robust estimation methods for regression modeling \cite{Jurecetal19,Alfoetal17}, model selection in robust regression is neglected in the literature. Model selection criteria  which are based on Gaussian assumptions discussed above are not robust since the log likelihood of normal distribution is used as lack of fit term.  First studies about the robustifed version of AIC were proposed by \cite{Hampel83model} and \cite{Ronchetti85}. For the penalty term, they use a generalization of AIC proposed by \cite{Bhandown1977} and given as the following form:
\begin{equation}\label{GAIC}
AIC=-2\log L(v;\hat{\boldsymbol{\beta}},\hat{\sigma}) + \kappa p.
\end{equation}
A choice for the parameter $ \kappa p$ in Eq. \eqref{GAIC} is given by \cite{Stone77} as follow:  
\begin{equation}\label{GAICL2L1}
AIC=-2\log L(v;\hat{\boldsymbol{\beta}},\hat{\sigma}) + \text{trace}(L_2^{-1}(\hat{\boldsymbol{\beta}},\hat{\sigma})L_1(\hat{\boldsymbol{\beta}},\hat{\sigma})),
\end{equation}
\noindent where 
\begin{equation}\label{L2type}
L_2(\boldsymbol \theta)=E\left(\frac{\partial^2}{\partial \boldsymbol \theta \partial \boldsymbol \theta^T}\log f(v;\boldsymbol \theta)\right),
\end{equation}
\begin{equation}\label{L1type}
L_1(\boldsymbol \theta)=E\left(\frac{\partial}{\partial \boldsymbol \theta}\log f(v;\boldsymbol \theta) \left(\frac{\partial}{\partial \boldsymbol \theta}\log f(v;\boldsymbol \theta) \right)^T \right).
\end{equation} 
$L_2$ and $L_1$ are also known as the definitions of Fisher Information Matrix (FIM). Let us introduce robust versions of Eqs. \eqref{L2type} and \eqref{L1type} using the objective function $\rho$ instead of  $\log(f)$, as considered by \cite{Hampel83model} and \cite{Ronchetti85}. 

\begin{equation}\label{L2typeR}
L_2(\boldsymbol \theta)=E\left(\frac{\partial^2}{\partial \boldsymbol \theta \partial \boldsymbol \theta^T}\rho(v;\boldsymbol \theta)\right),
\end{equation}
\begin{equation}\label{L1typeR}
L_1(\boldsymbol \theta)=E\left(\frac{\partial}{\partial \boldsymbol \theta}\rho(v;\boldsymbol \theta) \left(\frac{\partial}{\partial \boldsymbol \theta}\rho(v;\boldsymbol \theta) \right)^T \right).
\end{equation} 
Due to fact that $\psi(v;\boldsymbol{\theta})=\frac{\partial}{\partial \boldsymbol{\theta}} \rho(v;\boldsymbol{\theta})$, we can rewrite Eqs. \eqref{L2typeR} and \eqref{L1typeR} as follows:
\begin{equation}\label{L2typeRP}
L_2(\boldsymbol \theta)=E\left(\frac{\partial}{\partial \boldsymbol \theta }\psi(v;\boldsymbol \theta)\right),
\end{equation}
\begin{equation}\label{L1typeRP}
L_1(\boldsymbol \theta)=E(\psi(v;\boldsymbol{\theta}) \psi(v;\boldsymbol{\theta})^T).
\end{equation} 
As a consequent, the robust penalty term for AIC is proposed by \cite{Hampel83model}:
\begin{equation}\footnotesize \label{Hampel83penalty}
\alpha_H=\alpha_R+[E(\frac{\partial}{\partial \boldsymbol{\theta}} \psi(v;\boldsymbol{\theta}))]^{-1}E(\psi(v;\boldsymbol{\theta}) \psi(v;\boldsymbol{\theta})^T)[E(\frac{\partial}{\partial \boldsymbol{\theta}} \psi(v;\boldsymbol{\theta}))]^{-1}p.
\end{equation}
and alternative proposition, $\alpha_R$, is introduced by \cite{Ronchetti85}
\begin{equation}\label{Ronc84penalty}
\alpha_R=[E(\frac{\partial}{\partial \boldsymbol{\theta}} \psi(v;\boldsymbol{\theta}))]^{-1}E(\psi(v;\boldsymbol{\theta}) \psi(v;\boldsymbol{\theta})^T)p.
\end{equation}

For the lack of fit part, they use  objective function which defines robust estimators instead of using the likelihood function in AIC to reach the robustness. Therefore, $AIC_{\text{Hampel}}$ (AIC$_H$) and $AIC_{\text{Ronchetti}}$ (AIC$_R$) can be defined as the following forms:
\begin{equation}\label{AICHamp}
AIC_{\text{H}}=2\sum_{i=1}^n \rho(v;\boldsymbol{\theta}) + \alpha_H,
\end{equation} 
\begin{equation}\label{AICRonc}
AIC_{\text{R}}=2\sum_{i=1}^n \rho(v;\boldsymbol{\theta}) + \alpha_R.
\end{equation} 
When the objective function $\rho$ is chosen as Huber's $\rho_H$, the derived forms of Eqs.  \eqref{Hampel83penalty} and \eqref{Ronc84penalty} are given in Appendix A. 

The robust version of  ICOMP called as RICOMP(IFIM) defined by \cite{YanLiu} is given as follow 
\begin{equation}
RICOMP(IFIM)=-2 \log L(\hat{\boldsymbol{\theta}_R}) + 2C_1 (\mathcal{F}_R^{-1}),
\end{equation}
\noindent where $\hat{\boldsymbol{\theta}}=(\hat{\boldsymbol{\beta}}_R,\hat{\sigma}_R^2)$  represents the robust estimate of the parameter vector which can be obtained from the M-estimation and $\mathcal{\hat{F}}_R^{-1}$ is Inverse of Fisher Information Matrix (IFIM) is given by
\begin{equation}
\mathcal{\hat{F}}_R^{-1}=\widehat{Cov}(\hat{\boldsymbol{\beta}_R},\hat{\sigma}_R^2)=\begin{bmatrix}
\widehat{Cov}(\hat{\boldsymbol{\beta}_R}) & 0 \\
0 & 2\hat{\sigma}_R^2 \\
\end{bmatrix}.
\end{equation}

\noindent  So, $C_1(\mathcal{\hat{F}}_R^{-1})$ is computed by
\begin{equation}
C_1(\mathcal{\hat{F}}_R^{-1})=\frac{s}{2}\log \left[s^{-1} tr(\hat{\mathcal{F}}_R^{-1})\right]- \frac{1}{2} \log |\hat{\mathcal{F}}_R^{-1}|.
\end{equation}
For the robustness of lack of fit part in ICOMP,  \cite{YanLiu} used the Least Favorable distribution instead of normal distribution in $\log L$ and proposed $RICOMP(IFIM)$ as follow:
\begin{eqnarray}\label{leastfavo}
RICOMP(IFIM)&=&n \log(2 \pi)  + 2n \log (\sigma) -2 \log w(x_i) \\ \nonumber
&+& 2 \sum_{i=1}^{n} \rho \left(\frac{w(x_i)}{\hat{\sigma}} (y_i - x_i^T \hat{\boldsymbol{\beta}}) \right) + 2C_1 (\mathcal{\hat{F}}_R^{-1}). \\ \nonumber
\end{eqnarray}
For the special case of the function $w(x_i)=1$, Eq. \eqref{leastfavo} drops to Huber's $\rho$ function \cite{Huber81}.

Another definition of robust ICOMP called as RICOMP$_M$ from M-estimation method is defined by \cite{Guney18}
\begin{equation}\label{RICOMPay}
RICOMP_M= 2 \sum_{i=1}^{n} \rho(v_i;\hat{\boldsymbol{\beta}_M},\hat{\sigma_M}) +2C(\Sigma_M),
\end{equation}
\noindent where  $C$ is complexity part of ICOMP. In this  paper, we use $\rho$ as  Huber's $\rho_H$ function and $C$ as maximal complexity of Bozdogan's $C_1$.

\subsection{Information theoretic model selection criterion based on Shannon's relative entropy }\label{relativeentropyandcomplexity}
Let us consider a continuous d-dimensional distribution with joint density function $f(v_1,v_2,\cdots,v_d)$ and marginal density functions $f_j(v_j)$, $j=1,2,\cdots,d$.

Shannon's marginal entropy is given by
\begin{equation}\label{shannon}
H(X_j) = -E [\log f_j(v_j)] = -\int f_j(v_j) \log f(v_j) \ud v_j
\end{equation}
\noindent and joint or global entropy is given by \cite{Shannon48}
\begin{eqnarray}\label{HjointH}
H(v_1,v_2,\cdots,v_d) &=& - E [\log f(v_1,v_2,\cdots,v_d)] \\ \nonumber
&=& - \int \cdots \int f(v_1,v_2,\cdots,v_d) \log f(v_1,v_2,\cdots,v_d) \ud v_1 \ud v_2 \cdots \ud v_d. 
\end{eqnarray}
This is also called the Shannon's complexity \cite{Riss89}. Kullback-Leibler is used Eqs. \eqref{shannon} and \eqref{HjointH}  to define information divergence against independence. It is given by \cite{Kullback68}
\begin{equation}
KL(f_j(v_j) || f(v_1,v_2,\cdots,v_d))=\sum_{j=1}^{d} H(v_j) - H(v_1,v_2,\cdots,v_d).
\end{equation} 
This is also known as the mutual information. It was used to introduce information theoretic measure of covariance complexity by using Gaussian distribution by  \cite{Emden71}. The initial definition of covariance complexity is

\begin{equation}\label{VanEmden}
C_0(\sigma_{jj},\Sigma)=\frac{1}{2}\sum_{j=1}^{p} \log(\sigma_{jj})-\frac{1}{2}\log|\Sigma|,
\end{equation}
where $ \sigma_{jj} $= $ \sigma_{j}^2 $ is the j-th diagonal element of $ \Sigma $ and $p$ is the dimension of $ \Sigma $. Another definition of complexity which is obtained by maximization of $C_0$ according to set of all orthonormal transformations was proposed by \cite{Bozdogan1988}:  
\begin{equation}\label{HamparsumC1}
C_{1}(\Sigma)=\frac{p}{2}\log\left[\frac{tr(\Sigma)}{p} \right]-\frac{1}{2}\log|\Sigma|.
\end{equation}

Since the complexities in Eqs. \eqref{VanEmden} and \eqref{HamparsumC1} are based on the multivariate Gaussian  distributional assumption, we can say that they are not robust. Thus, we expect that they do not work properly in a Robust Model Selection Criteria (RMSC). In this paper, we introduce a new robust complexity term based on the structure of complexity in Eq. \eqref{VanEmden} by using Huber's $\rho_H$ objective function. Next section represents the theoretical details. 

\subsection{Huber's $\rho_H$ complexity based on KL divergence for the normal distribution}\label{hubercankayapamukcu}
Any objective function $\rho$ can be used  instead of using  $-\log(f)$ in Eq. \eqref{shannon} \cite{Huber81,phdcankaya15}. We use this similarity in construction  of our complexity term.  Let us redefine the Shannon's entropy based on $\rho$.
\begin{eqnarray}\label{Shannonrho}
H(v)&=&\int   -\log[f(v)] f(v) \ud v \\ \nonumber
H^{\rho}(v)&=&\int  \rho(v) f(v) \ud v,
\end{eqnarray} 
\noindent where $f$ is an underlying distribution. 

To obtain Huber's $\rho_H$ complexity ( $C_0^{\rho_H}$) based on Shannon's entropy, $\rho$ and $f$ in Eq. \eqref{Shannonrho} are chosen as $\rho_H$ and normal distribution, respectively.

Firstly, we consider to get univariate case, that is, $\sum_{j=1}^{p} E_g(\rho(u))=\int \rho(u) g(u) \ud u  $. It should be noted that the underlying distribution $g$ for expectation $E(\cdot)$ is chosen as standard normal distribution. Considering Huber's $\rho_H$ given as,
\begin{equation}\label{Huberrho}
\rho_H(u)=
\left\{
\begin{array}{ll}
u^2~ &,  |u| \leq k \\
2k|u|-k^2 &, |u| > k,
\end{array}
\right.
\end{equation}
\noindent we have
\begin{equation}\label{Huberintegral}
E(\rho_H(u)) = \int_{-\infty}^{-k} \left(-ku-\frac{k^2}{2} \right)g(u) \ud u + \int_{-k}^{k} u^2 g(u) \ud u + 
\int_{k}^{\infty} \left(ku-\frac{k^2}{2} \right) g(u) \ud u. 
\end{equation}
As is seen, there are three parts in Eq. \eqref{Huberintegral}. If we take,  $u=\frac{v-\mu}{\sigma}$. Then $\ud u=\ud v/\sigma$. Thus,
\begin{eqnarray}\label{Huberintegralms}
E(\rho_H(v;\mu,\sigma)) &=& \int_{-\infty}^{-k} \left(-k (\frac{v-\mu}{\sigma})-\frac{k^2}{2}\right)  f(v;\mu,\sigma)   \ud v  \\ \nonumber
&+& \int_{-k}^{k} \left(\frac{v-\mu}{\sigma}\right)^2 f(v;\mu,\sigma)   \ud v \\ \nonumber 
&+& \int_{k}^{\infty} \left(k(\frac{v-\mu}{\sigma})-\frac{k^2}{2}\right) f(v;\mu,\sigma)  \ud v, 
\end{eqnarray}
\noindent where $f(v;\mu,\sigma)= \frac{1}{\sigma\sqrt{2\pi}}\exp\{-\left(\frac{v-\mu}{\sqrt{2}\sigma}\right)^2\} $.

For calculating the integrals in Eq. \eqref{Huberintegralms}, we use the integral kernels in \cite{Canent18}. Then, we have
\begin{eqnarray}
E(\rho_H(v;\mu,\sigma)) &=& \frac{1}{\sigma \sqrt{\pi}}\left[\gamma(\frac{3}{2},\frac{k^2}{2})-\frac{k^2}{2}
\Gamma(\frac{1}{2},\frac{k^2}{2}) \right],
\end{eqnarray}
\noindent where $\gamma$ is lower and $\Gamma$ is upper incomplete gamma functions, respectively. $k$ is the tuning parameter of $\rho_H$.

Let us consider the multivariate case, we can write
\begin{equation}\label{multiintrho}
E(\rho_H(v_1,v_2,\cdots,v_n))=\int \rho_H(\left((v-\boldsymbol{\mu})^T\Sigma^{-1}(v-\boldsymbol{\mu})\right)^{1/2}) f(v;\boldsymbol{\mu},\Sigma) \ud v,
\end{equation}
where the multivariate form  $\left((v-\boldsymbol{\mu})^T\Sigma^{-1}(v-\boldsymbol{\mu})\right)^{1/2}$ can be  equivalent to $u=\frac{v-\mu}{\sigma}$ in the univariate case. 

For the integral of middle part in Eq. \eqref{multiintrho}, we have to calculate
\begin{equation}\label{multiintmiddle}
\int \frac{1}{2}(v-\boldsymbol{\mu})^T\Sigma^{-1}(v- \boldsymbol{\mu})f(v; \boldsymbol{\mu},\Sigma) \ud v,
\end{equation} 
where $  f(v; \boldsymbol{\mu},\Sigma)=(2\pi)^{-p/2}|\Sigma|^{-1/2}\exp\{-\frac{1}{2}[(v-\boldsymbol{\mu})^T\Sigma^{-1}(v-\boldsymbol{\mu})] \} $.

Let $u^2$ be $(v-\boldsymbol{\mu})^T \Sigma^{-1} (v-\boldsymbol{\mu})$. After taking the derivatives, $u \ud u=(v-\boldsymbol{\mu})\Sigma^{-1} \ud v$ and $u^{-2} (v-\boldsymbol{\mu})^T \Sigma^{-1} u \ud u = \ud v $. With little rearrangement, we have
\begin{equation}\label{starstep}
\Sigma (v-\boldsymbol{\mu})^{-1}u \ud u=\ud v.
\end{equation} 
If we take $v - \boldsymbol{\mu}=z$, we can write $u^2=z^T \Sigma^{-1} z$. Then $u \ud u=z \Sigma^{-1} \ud z$. Thus, we can rewrite Eq. \eqref{starstep} as follow,
\begin{equation}\label{starstep2}
\Sigma z^{-1}z \Sigma^{-1} \ud z = \ud v.
\end{equation} 
By using the integral kernels in Ref. \cite{Canent18}, we obtain following result for the integral of middle part in Eq. \eqref{multiintrho} given by
\begin{equation}\label{C0Hubermultimiddleresult}
\frac{1}{|\Sigma|^{1/2}(2\pi)^{p/2}} \sqrt{2}\gamma(\frac{3}{2},\frac{k^2}{2}).
\end{equation} 
For the tail part of $\rho_H$ function, we have to calculate following integral.
\begin{equation}\label{multiinttail}
\int (k|(v-\boldsymbol{\mu})^T\Sigma^{-1}(v-\boldsymbol{\mu})|^{1/2} - \frac{k^2}{2})  f(v; \boldsymbol{\mu},\Sigma) \ud v.
\end{equation} 
The absolute term $|(v-\boldsymbol{\mu})^T\Sigma^{-1}(v-\boldsymbol{\mu})|^{1/2}$ will produce the negative value if $v < \boldsymbol{\mu}$. Otherwise, it will be positive. So the result of left of integral is
\begin{equation}\label{tailabs}
\int k |(v-\boldsymbol{\mu})^T\Sigma^{-1}(v-\boldsymbol{\mu})|^{1/2} g(v; \boldsymbol{\mu},\Sigma) \ud v=0,
\end{equation}
\noindent and right part of integral is
\begin{equation}\label{tailk2}
-\int \frac{k^{2} }{2} g(v; \boldsymbol{\mu},\Sigma) \ud v =-\frac{1}{|\Sigma|^{1/2}(2\pi)^{p/2}} \frac{k^2}{\sqrt{2}}\Gamma(\frac{1}{2},\frac{k^2}{2}). 
\end{equation}

Note that the variable transformation also includes to change the limit values of integrals. So we omit to write the limits of integrals. 

Finally, we define our new robust penalty term as follow:
\begin{eqnarray}\label{C0Huber}
C_{0}^{\rho_H}(\Sigma(\hat{\boldsymbol \beta}),k)&=&\sum_{j=1}^{d} \frac{1}{\sigma_{jj}\sqrt{\pi}}\left[\gamma(\frac{3}{2},\frac{k^2}{2})-\frac{k^2}{2}
\Gamma(\frac{1}{2},\frac{k^2}{2}) \right] \\ \nonumber
&&-\frac{1}{|\Sigma|^{1/2}(2\pi)^{p/2}} \left[\sqrt{2}\gamma(\frac{3}{2},\frac{k^2}{2})-\frac{k^2}{\sqrt{2}}\Gamma(\frac{1}{2},\frac{k^2}{2}) \right].
\end{eqnarray}
The first part in Eq. \eqref{C0Huber} is sum of the marginal Shannon's entropy based on $\rho_H$ and the second part is joint Shannon's entropy  based on $\rho_H$. 
\subsection{A new proposition for the values of Huber's tuning parameter using $C_{0}^{\rho_H}$}\label{kc0new}
It is important that $ C_{0}^{\rho_H}$ depends on Huber's tuning parameter $k$. As is well-known, the value of tuning parameter $k$ is equal to $1.345$ in literature of robustness \cite{Andrewsetalloc72,Huber81}.  However, for $\Sigma=I_{(p \times p)}$,  if we take $k=1.345$, we obtain  $C_{0}^{\rho_H} \neq 0$ ($0.632784$). In this case is against to the definition of complexity of a covariance matrix. Additionally, since $C_{0}^{\rho_H}$ depends on the incomplete gamma functions, it does not equal to zero exactly. So we have to find which value of $k$ will give us $C_{0}^{\rho_H}(I_{(p \times p)}) \approx 0$. By using the simple grid search, we find $k=0.8875916$  Finally, $C_{0}^{\rho_H}$ has following properties:
\begin{itemize}
	\item $C_{0}^{\rho_H}(I_{(p \times p)}) \approx 0$, if $\Sigma=I_{(p \times p)}$ and $k=0.8875916$.
	\item $C_{0}^{\rho_H}(\Sigma) \rightarrow \infty$, if $|\Sigma| \rightarrow  0$ for all values of $k$.
\end{itemize}

This is an important finding and we call this new value of $k$ as  $k_{C_0}$.  It will be discussed how to obtain good estimates of parameters $\boldsymbol{\beta}$ using $k_{C_0}$  in the simulation section.

\subsection{A new robust information complexity criterion: $RICOMP_{C_0^{\rho_H}}$}\label{ourICOMP}
In this section, we define our new robust information complexity criterion ($RICOMP_{C_0^{\rho_H}}$). The robust penalty term, that is the second component of  $RICOMP_{C_0^{\rho_H}}$, is obtained  in Section \ref{hubercankayapamukcu}. For the robustness of first component of $RICOMP_{C_0^{\rho_H}}$, that is lack of fit term, we use the objective function instead of $-\log L$ in the lack of fit part by inspiring from \cite{Hampel83model,Ronchetti85,YanLiu}.  Therefore, our new criterion is defined as follow
\begin{equation}\label{ICOMPHuber}
RICOMP_{C_0^{\rho_H}}= n \log(2 \pi) + n \log(\hat{\sigma}^2) + 2 \sum_{i=1}^{n} \rho_H(v_i;\hat{\boldsymbol{\beta}},\hat{\sigma},k) + 2C_0^{\rho_H}(\Sigma).
\end{equation}
Eq. \eqref{ICOMPHuber} has four components. These are $ n \log(2 \pi)$, $n \log(\sigma^2) $, $\rho_H(v_i;\hat{\boldsymbol{\beta}},\hat{\sigma},k) $ and $C_0^{\rho_H}(\Sigma)$.  $ n \log(2 \pi)$, $n \log(\sigma^2) $ are from $RICOMP$ in Eq. \eqref{leastfavo} constructed based on the Least Favorable distribution to add scale parameter into model selection.  $n \log(\sigma^2) $ is also a part for lack of fit. Another part for lack of fit as a main part is represented by $\rho_H(v_i;\hat{\boldsymbol{\beta}},\hat{\sigma},k)$ to perform robust modeling. The last component is robust penalty term based on the Huber's objective function. Note that since we use $\rho_H$ for lack of fit part, the penalty term should be based on $\rho_H$. 

\section{Numerical examples}\label{numericalexamples}

\subsection{Monte Carlo simulation study for $k_{C_0}$}\label{MCkC0}

We discuss how $k_{C_0}$ will affect the robust estimation of parameters in the robust regression. 

We show  the simulation protocol to generate the explanatory variables $x=(x_1,x_2,\cdots,x_5)$ which are dependent each other, as given by the following design: 
\begin{itemize}
	\item $x_i=\sqrt{1-\alpha^2} z_i+ \alpha z_6$, $i=1,2,\cdots,5$
\end{itemize}
\noindent where the variables $z_1,z_2,\cdots,z_6$ are distributed from Burr III distribution with $c=2$ and $k=20$ as shape parameters. Thus, the variables $x_1,x_2,\cdots x_5$ have big values due to the heavy-tailedness property of BIII distribution \cite{Canetal19}. In our case, $\alpha=0.9$.   

Let us introduce how to generate random variables including outliers for the dependent variable $y$ in the regression model.
\begin{itemize}
	\item We have two populations with the sample sizes $n_1$ and $n_2$ from $N_1(\mu_1=0,\sigma_1=1)$ and $N_2(\mu_2=5,\sigma_2=10)$. We also use other contamination as $N_2(0,50)$. We combine these artificial data sets, i.e., $n=(n_1,n_2)$ with different percent of the Level of Contamination (LC) such as $\%10,\%20,\%30$ and $\%50$, i.e., $(1-\tau)N_1(0,1)+\tau N_2(5,10)$ and $(1-\tau)N_1(0,1)+\tau N_2(0,50)$. $\tau \in [0,1]$ is LC. The population  from $N_2$ is regarded as contamination into the population from $N_1$. These two contaminations $N_2(\mu_2=5,\sigma_2=10)$ and $N_2(\mu_2=0,\sigma_2=50)$  are used to generate outliers which are distant from underlying distribution $N_1(\mu_1=0,\sigma_1=1)$. Thus, we test the robustness property of $RICOMP_{C_0^{\rho_H}}$ in Eq. \eqref{ICOMPHuber}. Thus, we have the random numbers for $\varepsilon$ with $n \times 1$. 
	\item We generate dependent random variable as a true model:  $y=x_1+x_2+x_3+ \sigma \varepsilon$, where $\sigma=1$.
	\item 'robustfit' module in MATLAB2015a$\copyright$ is used to estimate regression and scale parameters robustly. 
	\begin{itemize}
		\item $\hat{Y_1}=X  \hat{\boldsymbol \beta}_R$ for $k=1.345$
		\item $\hat{Y_2}=X  \hat{\boldsymbol \beta}_R$ for $k_{C_0}=0.8875916$.
		\item  In order to compare the prediction performance of $\hat{Y_1}$ and $\hat{Y_2}$, we use mean absolute error (MAE), given by
		\begin{equation}
		\text{MAE}_j= \frac{1}{r} \sum_{t=1}^{r} \left( \frac{1}{n} \sum_{i=1}^{n}   |Y_i - \hat{Y}_{j_i}| \right) ~~ \text{for}~~ j=1,2.
		\end{equation}
	\end{itemize}
	where $r$ and $n$ are the numbers of replication and sample size, respectively.  Using MAE for testing the performance of prediction is better than  mean squared error, i.e., $ \frac{1}{r} \sum_{t=1}^{r} \left( \frac{1}{n} \sum_{i=1}^{n}   (Y_i - \hat{Y}_{j_i})^2 \right) $ \cite{Shao03}.  
\end{itemize} 
\begin{table}[!htb]\label{MAEtables}\centering
	\caption{Performance comparison of $k$ and $k_{C_0}$ in robust regression} 
	\begin{tabular}{cccccccc} \hline
		& 	& 	& &   \multicolumn{2}{c}{\tiny $(1-\tau)N_1(0,1)+\tau N_2(5,10)$}   & \multicolumn{2}{c}{\tiny $(1-\tau)N_1(0,1)+\tau N_2(0,50)$}   \\  
		n& LC* ($\%$)	& $n_1$	& $n_2$ & MAE$_1$  &  MAE$_2$ & MAE$_1$ & MAE$_2$   \\ \hline
		\multirow{4}{*}{n=30}
		& 10& 27 &3  & 11.020 & 10.861& 13.707&13.521 \\
		& 20 & 24 & 6 &11.334  &11.186 &16.802 &16.572 \\
		&30  & 21 & 9 &11.673  &11.518 &19.794 &19.504 \\
		&50  & 15 &15  & 12.315 &12.184  &25.595 &25.277\\   \hline
		\multirow{4}{*}{n=50} 
		&10  &45  & 5 &11.312  &11.165 &14.111 &13.950 \\
		&20  &40  &10  &11.822  & 11.674&17.227 & 17.038\\
		&30  & 35 & 15 &12.113  &11.977 &20.379 &20.166  \\
		&50  & 25 & 25 &12.912  &12.793 &26.507 &26.270 \\		\hline
		\multirow{4}{*}{n=100} 
		&10  & 90 &10  &11.653  &11.512 &14.359 & 14.207 \\
		&20  & 80 &20  &12.034  &11.396 &17.675 & 17.513\\
		&30  &  70& 30 &12.474  &12.340 &20.799 &20.627\\
		&50  &50  &50  & 13.274 &13.163 &27.131 &26.961 \\  \hline
	\end{tabular} \\
	\begin{tabular}{@{}c@{}}
		\multicolumn{1}{p{\textwidth}}{\footnotesize *LC represents Level of Contamination. MAE$_1$ and MAE$_2$ represent MAE of $k=1.345$ and $k_{C_0}=0.8875916$, respectively. $\tau \in [0,1]$ is LC.} 
	\end{tabular}
\end{table}
All computations were done on a MacBook Pro with 8 GB Memory, 3.8 GHz, IntelCore 8th generation i5. The results of MAE$_j$ are given in Table 1 
The number of replication is $r=10000$. According to Table 1, 
it is clear that robust regression models with $k_{C_0}$ has superior performance by yielding  less MAE. We propose that $k_{C_0}$ can be used as a new tuning parameter in robust modeling. Therefore, we will use  $k_{C_0}$ in the computation of robust regression models in the next simulation in order to be able to compare the model selection criteria.

\subsection{Monte Carlo simulation study for $RICOMP_{C_0^{\rho_H}}$}
In this section, we present Monte Carlo simulation study to show the performance  of $RICOMP_{C_0^{\rho_H}}$ in choosing true model by comparing other robust information criteria. The simulation protocol for generating the explanatory and dependent variables is same with previous simulation. We only take $N_2(\mu_2=5,\sigma_2=10)$ for the population with sample size $n_2$ as outliers. 

After having the artificial data sets, we use following steps in order to select the best model in competing models using RMSC.

\begin{itemize}
	\item All subset models constructed from the explanatory variables $x_1,x_2,\cdots,x_5$ are used. There are $2^{5}-1=31$ models to compare in competing subset models and we know that only subset of $\{x_1,x_2,x_3\}$ is true model.
	\item Robust regression are performed to these competing models. 
	\item  For each regression models, all RMSC are computed and listed. 
	\item The smallest value in the list is chosen.
	\item The models which correspond to the smallest model selection values are checked whether or not they are equivalent to true model.
	\item If they are the true model, then the criteria which select true model have earned one point and this score is stored. 
\end{itemize} 

\noindent	This process for all scenarios with different LC levels is repeated 10000 times. The results are given in Table 2. 

\begin{table}[!htb]\label{modelselecttables}\centering \tiny
	\caption{The number of choosing true model for the model selection criteria in 10000 runs} 
	\begin{tabular}{ccccccccc} \hline
		& LC* ($\%$)	& $n_1$	& $n_2$ & $RICOMP_{C_0^{\rho_H}}$  &   $RICOMP(IFIM)$ & $RICOMP_{M}$ & $AIC_{\text{H}}$ &$AIC_{\text{R}}$ \\ \hline
		\multirow{4}{*}{n=30}
		& 10& 27 &3  & 2798 & \bf{2801} &273 &1262 & 195\\
		& 20 & 24 & 6 &\bf{2637}  & 2117 &214 &1101 &146 \\
		&30  & 21 & 9 & \bf{2885} & 1946 & 147& 1395&79 \\
		&50  & 15 &15  &  \bf{2281}& 1604   &107 &1388 &24 \\   \hline
		\multirow{4}{*}{n=50} 
		&10  &45  & 5 &3082 & \bf{3505} &158 &564 &115 \\
		&20  &40  &10  & \bf{3108}  &2501   &100 &397 &77 \\
		&30  & 35 & 15 & \bf{3076} &2050 &62 &392 & 51 \\
		&50  & 25 & 25 & \bf{3665} & 2362  &31 & 971&17 \\		\hline
		\multirow{4}{*}{n=100} 
		&10  & 90 &10  &  3243& \bf{4347} &42 &134 &37  \\
		&20  & 80 &20  & \bf{3242} &  2966  & 6 &44 & 7\\
		&30  &  70& 30 & \bf{3198} &2146  &8 &45 &4 \\
		&50  &50  &50 & \bf{4535} & 3088   &1 &89 & 2\\  \hline
	\end{tabular} \\
	\begin{tiny}
		\begin{tabular}{@{}c@{}}
			\multicolumn{1}{p{\textwidth}}{\footnotesize *LC represents Level of Contamination. $AIC_{\text{H}}: AIC_{\text{Hampel}}$, $AIC_{\text{R}}$: $AIC_{\text{Ronchetti}}$ } 
		\end{tabular}
	\end{tiny}
\end{table}
\begin{itemize}
	\item The bold number represents the model selection criterion which has the highest number of choosing   the true model.
	\item For all scenario, $RICOMP_{C_0^{\rho_H}}$ has the highest number of choosing true model for the different LC except for $\%10$. It is observed that $RICOMP_{C_0^{\rho_H}}$ and $RICOMP(IFIM)$ have close performance in choosing the true model for LC$=\%10$.  Due to fact that LC is low,  it can be noted that the empirical distribution of data set is near to normal. Since the penalty term of  $RICOMP(IFIM)$ is based on the normal distribution, this result is expected.  
	\item Evaluating the result of $AIC_{\text{H}}$ in itself, it works well for $n=30$. However, its performance decreases noticeably when the number of observations increases.
	\item It is generally observed that $RICOMP_{C_0^{\rho_H}}$ and $RICOMP(IFIM)$ show better performance when compared other criteria.
	\item  $RICOMP_{C_0^{\rho_H}}$  yields more stable result even if LC increases.  It can be said that  $C_0^{\rho_H}$ has important role in robust model selection when the  outliers exist in data sets. 
	\item  For all scenario, $RICOMP_{M}$, $AIC_{\text{H}}$ and $AIC_{\text{R}}$ have poor  performance when compared others.
\end{itemize}
\subsection{Real data example}\label{realdata}
We present the performance of $RICOMP_{C_0^{\rho_H}}$ and the other model selection criteria to select the explanatory variables which can affect to the  response variable. For this purpose, we use the bridge construction data in \cite{Roozbehbuilt}. The variables in data are as follows:
\begin{itemize}
	\item Time: Time of construction ($y$)
	\item CCost: Cost of construction ($x_1$)
	\item Dwgs: Structural drawings number ($x_2$)
	\item Length: Length of bridge ($x_3$)
	\item Spans: Spans number ($x_4$)
	\item DArea: Deck of area of bridge ($x_5$)
\end{itemize}
The investigation about relationship among variables of data set was performed by \cite{Roozbehbuilt} and  illustrated  the scatter plot of data.  \cite{Roozbehbuilt} noted that there exists the striking nonlinear pattern among the variables, some of variables seem to be skewed and there are multicollinearity among some explanatory variables. After they applied log transformation on data set, they estimated regression parameters by using the Modified Least Trimmed Square Counter Multicollinearity (MLTSCM) as a robust regression method.   From this point, the regression techniques used by us and them for estimating regression parameters  are similar. These techniques can be comparable in the framework of efficient and robust modeling on the dependent variable $y$. However, applying the logarithmic transformation on the variables of a real data set may change the information. So it may not a good way to regularize the variables for modeling. If it is possible, modeling by using raw data without transformation on the variables may yield better result because type of transformation is not always one-to-one. In our case, we do not perform transformation and we analyze the data for all subset competing models by using robust regression with Huber M-estimator where $k=k_{C_0}=0.8875916$. For all competing models, RMSC are computed and assigned the best model according to RMSC's smallest values. The results are in  Table 3. 

\begin{table}[!htb]\small 
	\label{RMSCMAE}
	\caption{The best subset of predictors of regression model using RMSC for the bridge construction data set.}
	\begin{tabular}{cccc}
		Variables	& RMSC & Regression Model & MAE \\ \hline
		$\{x_1\}$	& $AIC_H:$ 3125.4 &$y=65.5715+0.2658 x_1$  &  48.9173 \\
		$\{x_4\}$  &$AIC_R:$ 31.205  & \multirow{2}{*}{$y=37.5115+45.8683x_4$}  &  \multirow{2}{*}{46.164 } \\
		$\{x_4\}$		&$RICOMP_M:$23.1171  &  & \\
		$\{x_1,x_4\}$	& $RICOMP(IFIM):$ 439.691  & $y=45.0571+0.1664x_1+22.7530x_4$ &44.9872  \\
		$\{x_2,x_5\}$	&$RICOMP_{C_0^{\rho_H}}:$ 415.2527 &$y=-13.0449+17.7851x_2+2.7736x_5$  &36.5115  \\
		\multirow{2}{*}{*$\{x_1,x_2,x_4\}$}	&  \multirow{2}{*}{-} & $\log (y)=2.0304 + 0.3056 \log(x_1) +$  & \multirow{2}{*}{36.146 }  \\ 
		&   & $ 0.6210 \log(x_2)+ 0.0657 \log(x_4)$  &   \\ \hline
	\end{tabular}
	\begin{tabular}{@{}c@{}}
		\multicolumn{1}{p{\textwidth}}{\footnotesize RMSC: Robust Model Selection Criteria. *Ref. \cite{Roozbehbuilt} } 
	\end{tabular}
\end{table}
According to Table 3, 
it is observed the best subsets of predictors are $S_1=\{x_2,x_5\}$ and $S_2=\{x_1,x_2,x_4\}$. Note that the subset of $S_1$ is obtained by using $RICOMP_{C_0^{\rho_H}}$ and the subset of  $S_2$ is provided by Ref. \cite{Roozbehbuilt}.

As it is seen, MAE values are very close for $S_1$ and $S_2$. However, it is important finding to obtain similar performance with fewer variables.  This shows utility and effectiveness of our proposed method with $RICOMP_{C_0^{\rho_H}}$.

Let us give a few comment about  performances of other RMSC. $RICOMP(IFIM)$ can be considered to have better performance than others. The worst performance is observed in $AIC_H$. We consider that it is due to penalty term in which has a part $x^Tx$ (see Appendix \ref{L2L1sectionproof}). As is well-known, when the multicollinearity problem  occurs in the data,  $x^Tx$ could be singular.

\begin{figure}[!htb]
	\centering 
	\begin{subfigure}{.48\linewidth}
			\centering
	\includegraphics[width=0.6\textwidth]{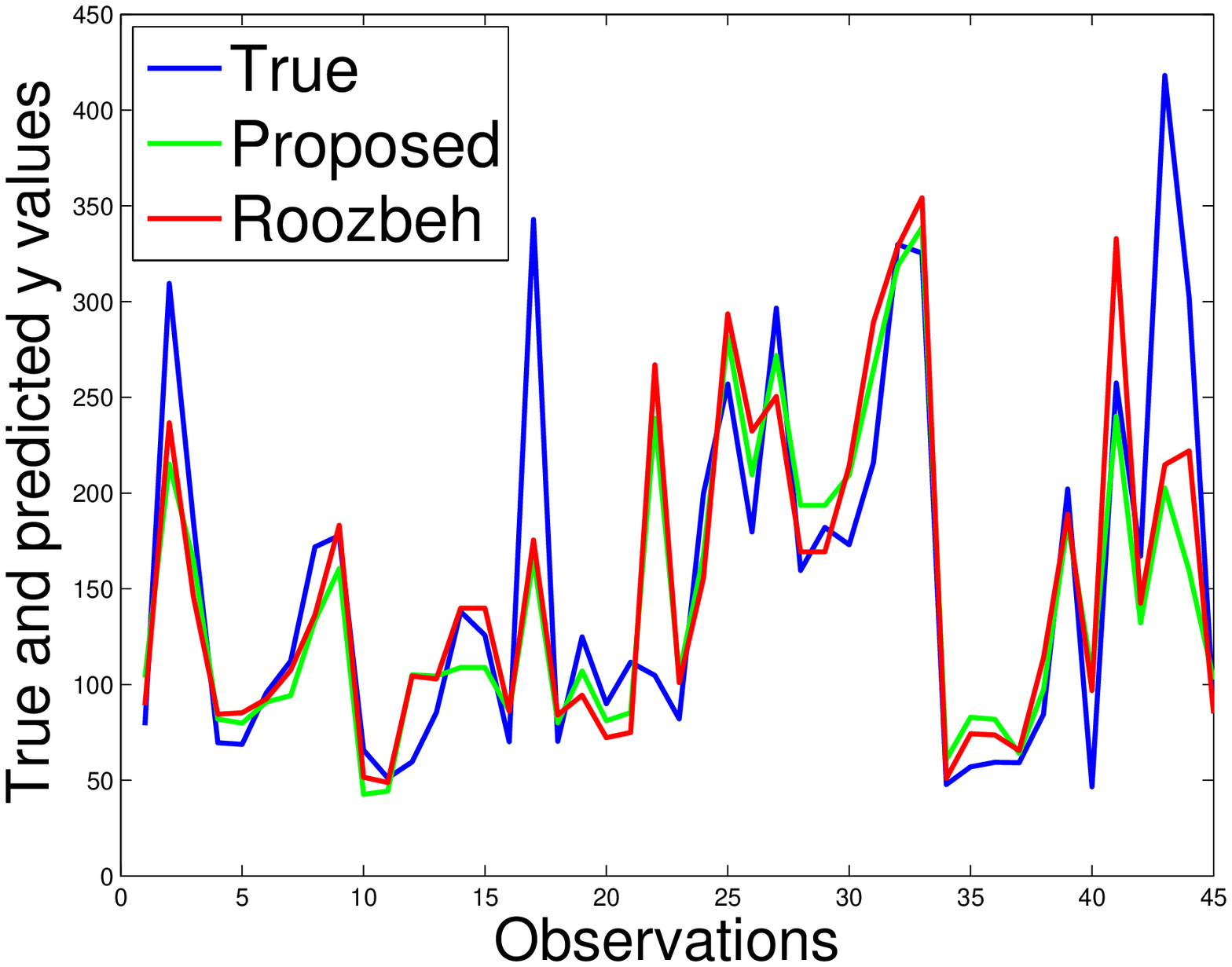}
		\caption{The true and predicted values for $y$}
	\end{subfigure}
	\begin{subfigure}{.48\linewidth}
			\centering
		\includegraphics[width=0.6\textwidth]{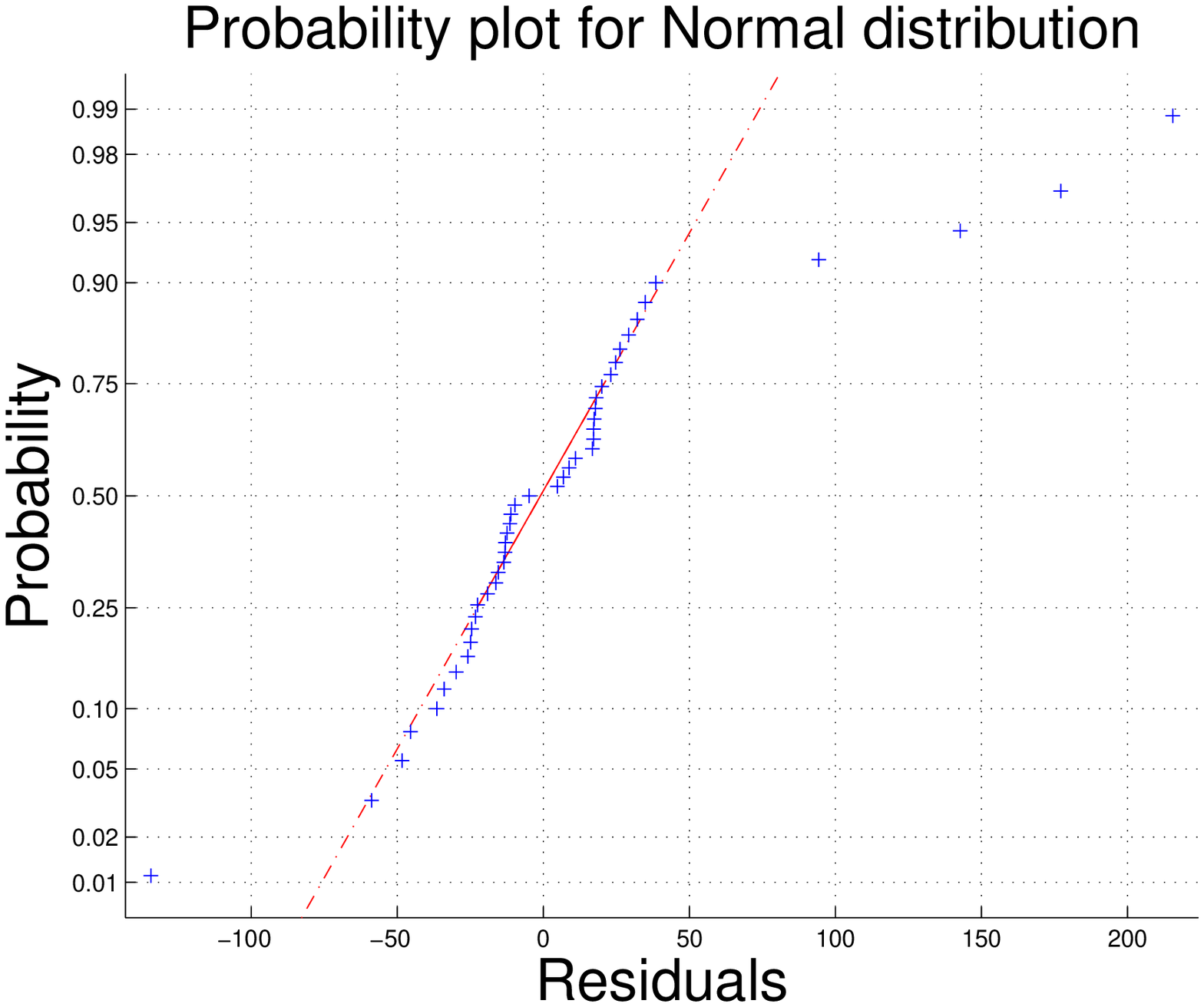}
		\caption{Probability plot of residuals from Huber's robust regression}
	\end{subfigure}
	\begin{subfigure}{.48\linewidth}
			\centering
	\includegraphics[width=0.6\textwidth]{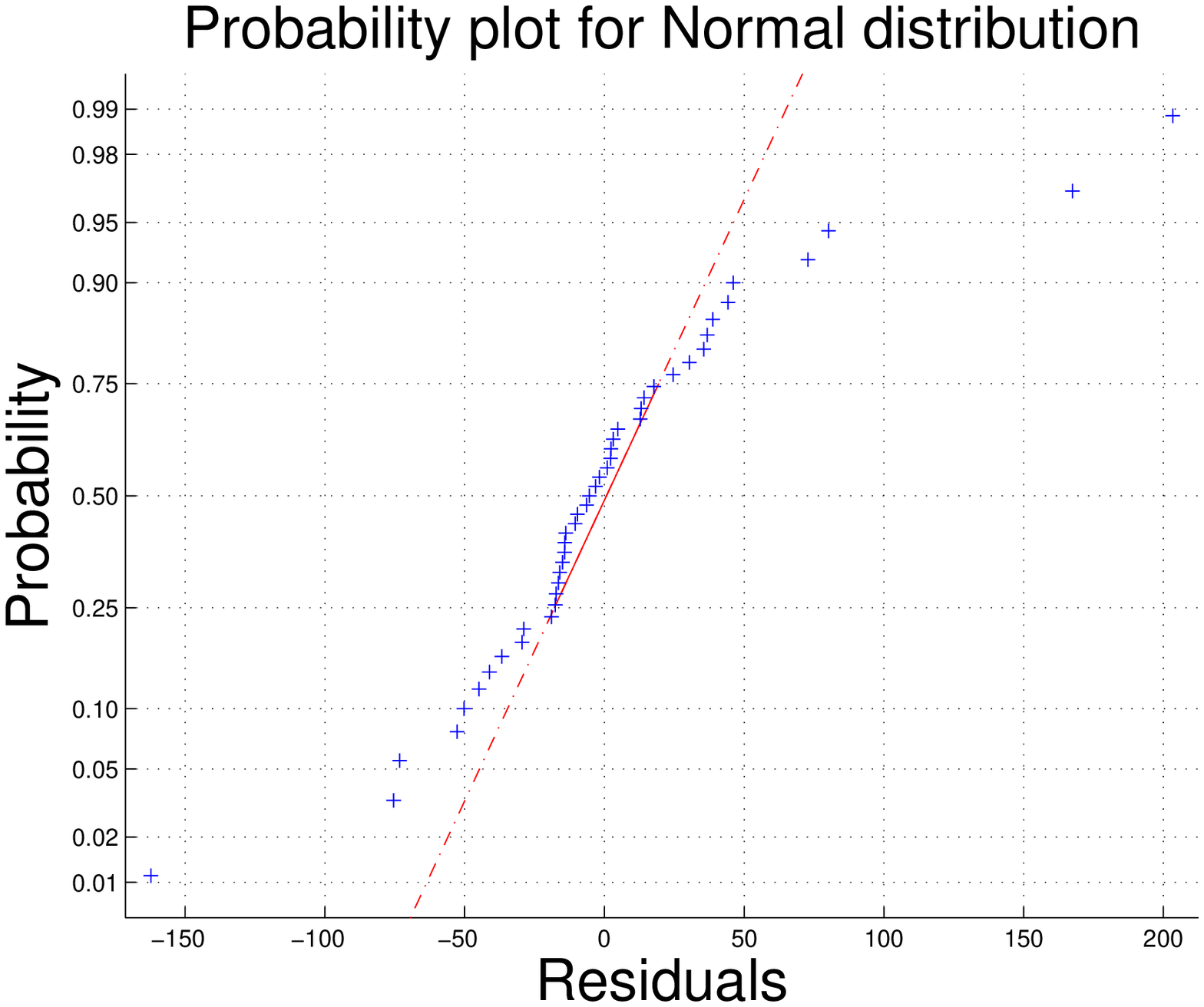}
		\caption{Probability plot of residuals obtained by  Ref. \cite{Roozbehbuilt}}
	\end{subfigure}
	\caption{(Color online) True and predicted values ($\hat{Y}$) from regression models  and probability plots of residuals}
	\label{Ex2figtoplucaAB}
\end{figure}


Fig. 1 (a) shows three lines to observe the nearness of prediction values obtained from regression models in last two lines in Table 3 to true values of $y$. Note that the values of MAE in Table 3 also show the prediction performance of regression models. Fig. 1 (b) and (c) show the goodness of fit by means of the probability plot.  The dashed line represents the probability values of normal distribution as the theoretical values. The blue pluses in Fig 1 (b) and (c) represent the empirical probabilities. When we compare the pluses in Fig. 1 (b) and (c), it is observed that the pluses in Fig. 1 (b) can fit the theoretical values better than that in Fig. 1 (c).

\section{Conclusion and discussion}\label{conclusion}
In this paper, we presented several forms of the robust model selection criteria. These are $AIC_H$, $AIC_R$, $RICOMP(IFIM)$ and $RICOMP_M$. We have extended information theoretic measure of complexity  ICOMP class of criteria for robust regression in order  to be resistant to high amount of outliers in data set. 

We used  Huber's $\rho_H$ to produce our robust penalty term for ICOMP and we introduced ${C_0^{\rho_H}}$. Then, we used as a specific form of the lack of fit part of $RICOMP(IFIM)$ by using Huber's objective function. Finally, we proposed a new robust information complexity criterion, $RICOMP_{C_0^{\rho_H}}$.
Additionaly, by using the properties of ${C_0^{\rho_H}}$, we proposed a new value of tuning parameter $k$ of Huber's $\rho_H$ as $0.8875916$.  We performed a Monte Carlo simulation to show important role of the value of $k$ in estimating regression parameters robustly. The simulation results  demonstrated clearly that our new value of $k$ called $k_{C_0}$ can be used as a tuning parameter in robust regression modeling. Moreover, $RICOMP_{C_0^{\rho_H}}$ was demonstrated on both Monte Carlo simulation study and real data to show the utility and effectiveness. According to numerical examples, the proposed criteria $RICOMP_{C_0^{\rho_H}}$ is successful to select the best model in robust regression for the data set which has high amount of outliers. In the future, $RICOMP_{C_1^{\rho_H}}$ will be derived from \cite{Bozdogan1988}.

 \acknowledgments{We would like to thank referees, Associative Editor(s) and Editor.  We would like to say our special thanks to Prof. Dr. Hamparsum Bozdogan for his constructive comments.}

\appendix

\section{Appendix: $F(\boldsymbol{\theta})$ and $R( \boldsymbol{\theta})$ of Huber's $\rho_H$ function }\label{L2L1sectionproof}
Eq. \eqref{GAICL2L1} includes the penalty term which is $\text{trace}(F^{-1}(\hat{\boldsymbol{\beta}},\hat{\sigma})R(\hat{\boldsymbol{\beta}},\hat{\sigma}))$.
We give the definition of $F$ and $R$ given by
\begin{equation}
F(\boldsymbol{\theta})=E\left(\frac{\partial^2 \rho(v;\boldsymbol{\theta})}{\partial \boldsymbol{\theta}   \partial \boldsymbol{\theta}^T} \right)=E\left(\frac{\partial}{\partial \boldsymbol{\theta}}\psi(v;\boldsymbol{\theta})\right)
\end{equation}
and 
\begin{equation}
R( \boldsymbol{\theta})=E \left( \left(\frac{\partial \rho(v;\boldsymbol{\theta})}{\partial \boldsymbol{\theta} }\right) \left(\frac{\partial \rho(v;\boldsymbol{\theta})}{\partial \boldsymbol{\theta} }\right)^T \right)=E\left(\psi(v;\boldsymbol{\theta})\psi(v;\boldsymbol{\theta})^T \right),
\end{equation}
\noindent where $ \boldsymbol{\theta}=(\boldsymbol{\beta},\sigma) $ \cite{Stone77,Ronchetti85}. In our case, we choose Huber's $\rho_H$ function for two cases. 

Firstly, let us give the derivatives of middle of Huber's $\rho_H$ with respect to (w.r.t.) parameters. 
\begin{equation}
\frac{\partial^2 \rho_H}{\partial \boldsymbol{\beta}  \partial \boldsymbol{\beta}^T}=\frac{x^Tx}{\sigma^2},
\end{equation}
\begin{equation}
\frac{\partial}{\partial \sigma}\frac{\partial \rho_H}{\partial \boldsymbol{\beta}}=\frac{\partial}{\partial \boldsymbol{\beta}}\frac{\partial \rho_H}{\partial \sigma}=\frac{2x}{\sigma^2}\left( \frac{y-x^T\boldsymbol{\beta}}{\sigma} \right),
\end{equation}
\begin{equation}
\frac{\partial^2 \rho_H}{\partial \sigma^2}=\frac{3}{\sigma^4} (y-x^T \boldsymbol{\beta})^2.
\end{equation}
Let us give the derivatives of tail part of Huber's $\rho_H$ w.r.t. parameters. 
\begin{equation}
\frac{\partial^2 \rho_H}{\partial \boldsymbol{\beta}  \partial \boldsymbol{\beta}^T}=0 x^Tx = \boldsymbol{0}_{p \times p},
\end{equation}
\begin{equation}
\frac{\partial}{\partial \sigma}\frac{\partial \rho_H}{\partial \boldsymbol{\beta}}=\frac{kx}{\sigma^2} \text{sign}(y-x^T \boldsymbol{\beta}),
\end{equation}
\begin{equation}
\frac{\partial^2 \rho_H}{\partial \sigma^2}=\frac{2k}{\sigma^3}|y-x^T \boldsymbol{\beta}|.
\end{equation}
Now we provide the $R$ case. For middle part, we have
\begin{equation}
\frac{\partial \rho_H}{\partial \boldsymbol{\beta}}=\frac{-x}{\sigma^2}(y-x^T \boldsymbol{\beta}),
\end{equation}
\begin{equation}
\frac{\partial \rho_H}{\partial \sigma}=\frac{-1}{\sigma^3}(y-x^T \boldsymbol{\beta})^2.
\end{equation}
We have products of the derivatives of $\rho_H$ w.r.t. parameters:
\begin{equation}
\frac{\partial \rho_H}{\partial \boldsymbol{\beta}} \frac{\partial \rho_H}{\partial \boldsymbol{\beta}}=\frac{x^T x}{\sigma^4}(y-x^T \boldsymbol{\beta})^2,
\end{equation} 
\begin{equation}
\frac{\partial \rho_H}{\partial \sigma} \frac{\partial \rho_H}{\partial \sigma}=\frac{1}{\sigma^6}(y-x^T \boldsymbol{\beta})^4,
\end{equation}
\begin{equation}
\frac{\partial \rho_H}{\partial \boldsymbol{\beta}} \frac{\partial \rho_H}{\partial \sigma}=\frac{x}{\sigma^5}(y-x^T \boldsymbol{\beta})^3. 
\end{equation}
Similarly, one can obtain the derivatives of tail part of Huber's $\rho_H$ function.

Secondly,  we give expectation for the elements of matrices. We use the prepared kernels from Ref. \cite{Canent18} to calculate the integrals of $F$ and $R$. The underlying distribution of expectation is chosen as normal distribution.
\begin{eqnarray}\label{L2bb}
E\left(\frac{\partial^2 \rho_H}{\partial \boldsymbol{\beta}  \partial \boldsymbol{\beta}^T} \right)&=& \int_{-\infty}^\infty \left(\frac{\partial^2 \rho_H}{\partial \boldsymbol{\beta}  \partial \boldsymbol{\beta}^T} \right) f(v;\boldsymbol{\beta},\sigma)   \\ \nonumber 
&=&\frac{x^Tx}{\sigma^2 \sqrt{\pi}}\gamma(1/2,k^2/2),
\end{eqnarray}
\noindent where $f(v;\boldsymbol{\beta},\sigma)= \frac{1}{\sigma\sqrt{2\pi}}\exp\{-\left(\frac{y-x^T\boldsymbol{\beta}}{\sqrt{2}\sigma}\right)^2\} $.
The undiagonal elements of matrix $F$ and $R$, that is,
\begin{equation}\label{L2bs}
E\left(\frac{\partial}{\partial \boldsymbol{\beta}}\frac{\partial \rho_H}{\partial \sigma }\right)=E\left(\frac{\partial}{\partial \sigma}\frac{\partial \rho_H}{\partial \boldsymbol{\beta} }\right)=E\left(\frac{\partial \rho_H}{\partial \boldsymbol{\beta}}\frac{\partial \rho_H}{\partial \sigma}\right)=0.
\end{equation}
Note that regression and scale parameters are independent for the regression and scale case in Eq. \eqref{L2bs}.  One can show these results easily.  
\begin{equation}\label{L2ss}
E \left(\frac{\partial^2 \rho_H}{\partial \sigma^2} \right)=\frac{6}{\sigma^2\sqrt{\pi}}\gamma(3/2,k^2/2)+\frac{2\sqrt{2}k}{\sigma^2\sqrt{\pi}}\Gamma(1,k^2/2).
\end{equation}
Eqs. \eqref{L2bb}-\eqref{L2ss} are for $F$ case. Now let us give the integral results for $R$ case.
\begin{equation}\label{L1bb}
E\left(\frac{\partial \rho_H}{\partial \boldsymbol{\beta}}\frac{\partial \rho_H}{\partial \boldsymbol{\beta}} \right)=\frac{x^Tx}{\sigma^2 \sqrt{\pi}}\left[2\gamma(3/2,k^2/2)+k^2\Gamma(1/2,k^2/2)\right],
\end{equation}
\begin{equation}\label{L1ss}
E \left(\frac{\partial \rho_H}{\partial \sigma}\frac{\partial \rho_H}{\partial \sigma}\right)=\frac{2}{\sigma^2\sqrt{\pi}}\left[2\gamma(5/2,k^2/2)+k^2\Gamma(3/2,k^2/2)\right].
\end{equation}

Eqs. \eqref{L2bb}-\eqref{L1ss} show the calculated form of Eq. \eqref{Hampel83penalty}. Let us give matrix elements of Eq. \eqref{Hampel83penalty}:
\begin{eqnarray}\label{H11}
A_{11}&=&\boldsymbol{1}_{p \times p} \left[2\gamma(3/2,k^2/2)+k^2\Gamma(1/2,k^2/2)\right]/\gamma(1/2,k^2/2) \\ \nonumber
&+&\frac{2\gamma(3/2,k^2/2)+k^2\Gamma(1/2,k^2/2)}{\gamma^2(1/2,k^2/2)}(x^Tx)^{-1}\sigma^2 \sqrt{\pi},
\end{eqnarray}
\begin{equation}\label{H12H21}
A_{12}=A_{21}=0,
\end{equation}
\begin{eqnarray}\label{H22}
A_{22}&=& \frac{2\gamma(5/2,k^2/2)+k^2\Gamma(3/2,k^2/2)}{3\gamma(3/2,k^2/2)+\sqrt{2}k\Gamma(1,k^2/2)} \\ \nonumber
&&\left(1+\frac{\sigma^2 \sqrt{\pi}}{2\left[3\gamma(3/2,k^2/2)+\sqrt{2}k\Gamma(1,k^2/2)\right]}\right). 
\end{eqnarray}
Eq. \eqref{H11} has $(x^Tx)^{-1}$. 
\begin{equation}
A=\begin{bmatrix}
A_{11} && A_{12} \\
A_{21} && A_{22} 
\end{bmatrix}.
\end{equation}

\begin{thebibliography}{}

\bibitem{Akaike73}	Akaike H (1973) Maximum likelihood identification of Gaussian autoregressive moving average models. Biometrika, 60(2), 255-265

\bibitem{AkbilgicBoz11} Akbilgic O, Bozdogan H (2011) Predictive subset selection using regression trees and RBF neural networks hybridized with the genetic algorithm. European Journal of Pure and Applied Mathematics, 4(4), 467-485

\bibitem{Alfoetal17}  Alf\`{o} M, Salvati N,  Ranallli MG (2017) Finite mixtures of quantile and M-quantile regression models. Statistics and Computing, 27(2), 547-570


\bibitem{Akaike74} 	Akaike H (1974) A new look at the statistical model identification. In Selected Papers of Hirotugu Akaike (pp. 215-222). Springer, New York, NY


\bibitem{Andrewsetalloc72}   Andrews DF,  Bickel  PJ,  Hampel FR,  Huber PJ,  Rogers WH,  Tukey JW (1972) Robust Estimates of Location: Survey and Advances, Princeton University, Princeton, N.J. 



\bibitem{Bhandown1977} Bhansali RJ,  Downham DY (1977) "Some properties of the order of an autoregressive model selected by a generalization of Akaike? s EPF criterion." Biometrika 64.3 : 547-551

\bibitem{Bozdogan87} Bozdogan H (1987) Model selection and Akaike's information criterion (AIC): The general theory and its analytical extensions. Psychometrika 52.3 : 345-370

\bibitem{Bozdogan1988} Bozdogan H (1988)  ICOMP: A new model selection criterion. In: Hans H. Bock (Ed.), Classification and Related Methods of Data Analysis, North-Holland, Amsterdam, April, 599-608

\bibitem{phdcankaya15} \c{C}ankaya, MN (2015) "M-Estimators with asymmetric influence function: Properties and their applications." PhD diss., University of Ankara

\bibitem{Canent18} \c{C}ankaya M (2018) Asymmetric bimodal exponential power distribution on the real line. Entropy, 20(1), 23

\bibitem{Canetal19} \c{C}ankaya M N, Yal\c{c}{\i}nkaya A, Alt{\i}nda\v{g} \"{O}, Arslan O (2019) On the robustness of an epsilon skew extension for Burr III distribution on the real line. Computational Statistics, 1-27

\bibitem{Emden71} Emden V (1971) An analysis of complexity. Mathematiscxh Centrum

\bibitem{Guney18} G\"{u}ney Y  (2018) Robust Model Selection Using Information Complexity (ICOMP) Crtiterion for Symmetric and Skew Distributions. PhD diss.,  Ankara University

\bibitem{Hampel83model} Hampel FR (1983) Some aspects of model choice in robust statistics. In Proceedings of the 44th Session of the ISI (Vol. 2, pp. 767-771

\bibitem{Hampeletal86}   Hampel FR, Ronchetti EM, Rousseeuw PJ, Stahel WA (1986) Robust Statistics: The Approach Based on Influence Functions, Wiley Series in Probability and Statistics, New York, USA



\bibitem{HanQuin1979} Hannan EJ,  Quinn BG (1979) The determination of the order of an autoregression. Journal of the Royal Statistical Society: Series B (Methodological) 41(2) 190-195


\bibitem{Huber81}  Huber PJ (1981) Robust Statistics, Wiley, New York

\bibitem{Howeetal12} Howe ED, Bozdogan H,  Kiroglu G (2012) Performance of information complexity criteria in structural equation models with applications. European Journal of Pure and Applied Mathematics, 5(3), 282-301

\bibitem{Jurecetal19}   Jure\v{c}kov\'{a} J, Picek J,  Schindler M (2019) Robust statistical methods with R. Chapman and Hall/CRC

\bibitem{Kullback68}  Kullback, S. (1968) Information Theory and Statistics. Dover, New York

\bibitem{KocBoz15} Koc EK, Bozdogan H. (2015) Model selection in multivariate adaptive regression splines (MARS) using information complexity as the fitness function. Machine Learning, 101(1-3), 35-58

\bibitem{Kwonetal08} Kwon Y, Bozdogan H, Bensmail H (2008) Performance of model selection criteria in Bayesian threshold VAR (TVAR) models. Econometric Reviews, 28(1-3), 83-101

\bibitem{Linzuc86} Linhart H, Zucchini W (1986) Model selection. John Wiley and Sons

\bibitem{YanLiu} Liu Y (2007) Robust and Misspecification Resistant Model Selection in Regression Models with Information Complexity and Genetic Algorithms

\bibitem{Pamukcu15} Pamukcu E, Bozdogan H, \c{C}al{\i}k S (2015) A novel hybrid dimension reduction technique for undersized high dimensional gene expression data sets using information complexity criterion for cancer classification. Computational and mathematical methods in medicine, 2015


\bibitem{Riss89} Rissanen J (1989) Stochastic complexity in statistical inquiry

\bibitem{Ronchetti85}	Ronchetti E (1985) Robust model selection in regression. Statistics \& Probability Letters 3 21-23

\bibitem{Schwarz78} Schwarz G (1978) Estimating the dimension of a model. The annals of statistics, 6(2), 461-464

\bibitem{Shannon48} Shannon CE (1948) A mathematical theory of communication. Bell Systems Technology Journal, 27, 379-423

\bibitem{Shao03} Shao J (2003)  Mathematical Statistics. Springer Texts in Statistics

\bibitem{Stone77}  Stone M (1977) An asymptotic equivalence of choice of model by cross-validation and Akaike's criterion. Journal of the Royal Statistical Society: Series B (Methodological), 39(1), 44-47

\bibitem{sugiura1978} Sugiura N (1978) Further analysts of the data by akaike's information criterion and the finite corrections: Further analysts of the data by akaike's. Communications in Statistics-Theory and Methods 7.1 13-26

\bibitem{takeuchi1976} Takeuchi K (1976)  The distribution of information statistics and the criterion of goodness of fit of models. Mathematical Science 153 12-18

\bibitem{Roozbehbuilt} Roozbeh M, Babaie-Kafaki S,  Sadigh, AN (2018) A heuristic approach to combat multicollinearity in least trimmed squares regression analysis. Applied Mathematical Modelling, 57, 105-120

\end{thebibliography}
\end{document}